\documentclass[12pt,reqno]{amsart}
\usepackage{amssymb}
\usepackage{bm}
\usepackage{graphicx}
\usepackage[centertags]{amsmath}
\usepackage{amsfonts}
\usepackage{amssymb}
\usepackage{amsthm}
\usepackage{enumerate}
\newtheorem{theorem}{Theorem}[section]
\newtheorem{definition}{Definition}[section]

\newtheorem{remark}{Remark}[section]

\newtheorem{lemma}[theorem]{Lemma}

\input epsf

\begin{document}

\title{Circle products as restrictions of the square product}

\author[V. A. Anagnostopoulos, Y. Sarantopoulos]{Vasileios A. Anagnostopoulos, Yannis Sarantopoulos}


\address[V. A. Anagnostopoulos]
{Department of Mathematics\\
National Technical University \\
Zografou Campus 157 80, Athens, Greece.}
\email{vanag@telecom.ntua.gr}

\address[Yannis Sarantopoulos]
{Department of Mathematics\\
National Technical University \\
Zografou Campus 157 80, Athens, Greece.}
\email{ysarant@math.ntua.gr }



\begin{abstract}
 We equip the tensor algebra of a vector space $U$ over the real or complex field with an alternative product. The new product has the property that if we specialize it to the symmetric tensor algebra becomes the circle product introduced by Brouder \cite{brouder09} while specialization to the antisymmetric algebra becomes the circle product introduced by Rota and Stein \cite{rota94}. We prove some interesting properties of this product analogous to the properties proved by the previous authors. We use Hopf algebraic methods to simplify our results. Finally we prove that the classical algebraic structure of the tensor algebra and the new one are homomorphically isomorphic.\bigskip
\end{abstract}

\maketitle

\section{Introduction}

During the last few years there has been a growing body of literature targeting the algebraic structure of quantum field theory. A purely algebraic interpretation of Wick ordering  \cite{fauser01} was provided giving insight to one of the most intriguing aspects of theoretical physics. On the other hand, this work provides more insight at the mathematical front to the structure of the tensor algebra of a vector space. Our interest lies on the mathematical basis of the approach and not its physical meaning.  Circle products were introduced by Rota and Stein in \cite{rota94} in the context of super-symmetric Hopf algebras. Brouder in \cite{brouder09} carried out the same construction to derive similar quantities for the symmetric algebra. The above work leaves out the important question of whether or not this construction is carried out for the symmetric tensor algebra in an analogous manner, like Rota's and Stein's construction, in the context of the antisymmetric algebra. Though the two constructions are similar in spirit there is no clear indication of an underlying more general construction. Since both algebras are quotient algebras of the tensor algebra the most natural generalization is to create a unified construction to the tensor algebra and then to specialize to subspaces. More specifically, one expects to see a circle product and a Laplace pairing on the tensor algebra that specializes to the above constructions when transferred to the quotient algebras. In this work we achieve exactly this. We carry out this construction on the tensor algebra and re-derive the two specialized constructions. In this interpretation, our work is complementary to both these works and unifies them. We re-derive the symmetric algebra as a necessary ingredient of this approach. We also re-derive the circle products of Brouder and Rota as consequences of more general constructions.

\section{Definitions, notation and basic results}
In this section we fix notation and review some basic facts on the Hopf algebraic structure of the tensor algebra of a finite dimensional vector space which will be used in the next sections. We also define the joint tensor algebra of two vector spaces in duality and describe its Hopf algebraic structure. Finally, we revisit the definition of a circle product defined on the symmetric algebra generated by a finite dimensional vector space equipped with an inner product. For basic definitions of a Hopf algebra, reference \cite{podles98} provides an excellent short introduction. See also \cite{majid02,kassel95} which are the standard textbooks on the subject among others. For the definition of the tensor algebra of a vector space and symmetric tensors we refer to \cite{cartier07}.

\subsection{Review of basic facts on Hopf algebras}
Let $A$ be a unital algebra over the field $\mathbb{K}$. In general its dual is not an algebra. If there is a special additional structure that we will outline below then not only its dual becomes a unital algebra but this additional structure is also present to the dual.
\begin{definition}
Let $A$ be a unital algebra with unit $\mathbf{1}$. Suppose that following mappings are defined
\begin{enumerate}
\renewcommand{\labelenumi}{(\roman{enumi})}
\item $\Delta : A \to A \otimes A$	(the co-multiplication, homomorphism),
\item $S : A \to A \otimes A$		(the antipode, anti-homomorphism),
\item $\epsilon : A \to \mathbb{K}$	(the co-unit, homomorphism)
\end{enumerate}
with the additional properties
\begin{enumerate}
\renewcommand{\labelenumi}{(\roman{enumi})}
\item $(I \otimes \Delta) \Delta= (\Delta \otimes I) \Delta$	(the co-associativity),
\item $(\epsilon \otimes \Delta) \Delta= (\Delta \otimes \epsilon) \Delta = I$,
\item  $(S \otimes I) \Delta= (I \otimes S) \Delta = \epsilon(\cdot) \mathbf{1}$	(the co-unit)
\end{enumerate}
where $I$ is the identity operator. We call the quadruple $\langle A, \Delta, S, \epsilon \rangle$  a \emph{Hopf algebra}.
\end{definition}

We use the well known Sweedler notation, (see \cite{kassel95}),  for the co-multiplication of an element $a$ of the Hopf algebra $A$ as
\begin{equation*}
\Delta(a)  =  \sum_{(a)} a_{(1)} \otimes a_{(2)}\\
\end{equation*}

We call a Hopf algebra co-commutative when
\begin{equation*}
\Delta(a)  =  \sum_{(a)} a_{(1)} \otimes a_{(2)} =  \sum_{(a)} a_{(2)} \otimes a_{(1)} \,.
\end{equation*}

For repeated application of $\Delta$, we use the standard notation,
\begin{equation*}
 (I\bigotimes \Delta) \Delta(a)=  \sum_{(a)} a_{(1)} \otimes a_{(2,1)} \otimes a_{(2,2)}\,
\end{equation*}
and
\begin{equation*}
 (\Delta\bigotimes I) \Delta(a)=  \sum_{(a)} a_{(1,1)} \otimes a_{(1,2)} \otimes a_{(2)}\,.
\end{equation*}

We call a Hopf algebra co-associative when the two quantities above are equal.

A standard result \cite{majid02} states that given two Hopf algebras $A_{1},A_{2}$ we can turn also $A_{1} \bigotimes A_{2}$ into a Hopf algebra. We summarize now the details. Since $A_{1},A_{2}$ are unital algebras their tensor product is again an algebra with product (defined for elementary tensors)
\begin{equation*}
(a_{1} \otimes a_{2} ) (a^{'}_{1} \otimes a^{'}_{2} )= a_{1} a^{'}_{1} \otimes a_{2}a^{'}_{2}\,.
\end{equation*}

Let us define for $a \in A_{1}$ and $b \in A_{2}$
\begin{enumerate}
\renewcommand{\labelenumi}{(\roman{enumi})}
\item $(\Delta_{1} \otimes \Delta_{2})(a_{1} \otimes a_{2} ) =  \sum_{(a),(b)} (a_{(1)i} \otimes b_{(1)}) \otimes (a_{(2)} \otimes b_{(2)})$\,,
\item $(S_{1} \otimes S_{2}) (a_{1} \otimes b )  =  S_{1}(a)S_{2}(b)$\,,
\item $(\epsilon_{1} \otimes \epsilon_{2}) (a \otimes b )  =  \epsilon_{1}(a) \epsilon_{2}(b)$\,,
\end{enumerate}

\begin{theorem}\label{thm:tensorhopf}
The quadruple $\langle A_{1}\bigotimes A_{2},\Delta_{1} \otimes \Delta_{2}, S_{1} \otimes S_{2},  \epsilon_{1} \otimes \epsilon_{2} \rangle$ is a Hopf algebra.
\end{theorem}
\subsection{Review of basic facts on tensor algebra}
 Let $U$ be a vector space. The letter $x$ possibly indexed by some subscript
will usually denote a generic element of U. Its tensor algebra $\mathcal{T}(U)$ is defined as the direct sum of the $n$-fold tensor products of $U$ with itself $\bigotimes^{n}U$, namely
\[
\mathcal{T}(U) = \bigoplus_{n=0}^{\infty}  \bigotimes ^{n} U\,.
\]
We define $\bigotimes^{0}U =\mathbb{K}$, the field of $U$. The unit element of the field when we view it as a member of the tensor algebra will be denoted by ${\bm 1}$.  For our purposes we restrict our attention to the complex field $\mathbb{C}$ or to the real field $\mathbb{R}$. One can also view the tensor algebra as the free non-commutative algebra generated by the  elements of $U$. Since $\bigotimes^{n}U$ is the span of elementary tensors $x_{1}\otimes \cdots \otimes x_{n}$ where $x_{i} \in U$, we define the product in the algebra between two elements of $\bigotimes^{n}U$ and $\bigotimes^{m}U$ as
\begin{equation*}
(x_{1}\otimes \cdots \otimes x_{n}) \cdot (x^{'}_{1}\otimes \cdots \otimes x^{'}_{m})=x_{1}\otimes \cdots \otimes x_{n} \otimes x^{'}_{1}\otimes \cdots \otimes x^{'}_{m}\,,
\end{equation*}
and extend with linearity to the whole algebra. It can be proved that with this multiplication the tensor algebra is a unital associative algebra. Because of the direct sum construction of the tensor algebra and the definition of the product, it is a graded algebra.

Let $U$ and $V$ be two vector spaces in duality via the bilinear form $\langle \cdot , \cdot \rangle$ , both over the same field $\mathbb{K}$. The letter $x$ possibly indexed by some subscript will usually denote a generic element of U while the letter $y$ indexed by some subscript will denote a generic element of V. Their {\it joint tensor algebra} is defined as the pair consisting of the duality and the tensor product
\[
\mathcal{T} (U, V)= \mathcal{T}(U) \bigotimes \mathcal{T}(U)\,.
\]
The term ``{\it joint tensor algebra}" is not used in the literature. However, since our work pertains to the properties of this specific tensor product we decided to give it a name encoding its special nature. It is easy to see that we can equip it with a suitable product to become a unital non-commutative associative algebra. We define it first for elementary tensors, that is
\begin{equation}
(u_{1} \otimes v_{1}) \cdot (u_{2} \otimes v_{2}) = (u_{1} \cdot u_{2}) \otimes (v_{1} \cdot v_{2})\,,
\end{equation}
where $u_{i} \in \bigotimes^{n_{i}}U$ and $v_{j} \in \bigotimes^{m_{j}}V$. We can extend it with linearity to the joint tensor algebra. It is well-known that the bilinear form can be extended to the joint tensor algebra.
The extension takes place by first defining it for elementary tensors
\begin{equation*}
\langle (x_{1}\otimes \cdots \otimes x_{n}) , (y_{1}\otimes \cdots \otimes y_{m})\rangle =\delta_{m,n} \langle x_{1}, y_{1} \rangle \cdots \langle x_{n}, y_{m} \rangle\,.
\end{equation*}
A typical element of the joint tensor algebra can be written as
\begin{equation*}
a= \sum_{i} u_{i}  \otimes v_{i}\,,
\end{equation*}
where $u_{i} \in \bigotimes^{n_{i}} U$ and $v_{i} \in \bigotimes^{m_{i}} V$.
For this reason we define the bilinear form for two elements in the joint tensor algebra $a,b$ as
\begin{equation*}
\langle a , b \rangle := \langle \sum_{i} u_{i}  \otimes v_{i} , \sum_{j} u^{'}_{j}  \otimes v^{'}_{j} \rangle=\sum_{i}  \sum_{j} \langle  u_{i} , v^{'}_{j} \rangle \langle  u^{'}_{j} , v_{i} \rangle\,.
\end{equation*}

We will specialize the previous definitions to an important quotient algebra of the tensor algebra, namely the symmetric tensor algebra. On the tensor algebra of the vector space $U$ we can define the symmetrization operator first for elementary tensors as
\begin{equation}\label{symm:def}
{\rm Symm}( x_{1} \otimes \cdots \otimes x_{n}) := \frac{1}{n!} \sum_{p \in Per_{n}}  x_{p(1)} \otimes \cdots \otimes \cdots x_{p(n)}
\end{equation}
and extend it by linearity on the whole tensor algebra. We collect in the next theorem some useful facts on the symmetrization operator \cite{rodrigues07,ryan80}.

\begin{theorem}\label{thm:symm}
Let $U$ be a vector space. Then, the symmetrization operator defined by (\ref{symm:def}) is a projection. The projection of the tensor algebra is denoted by $\mathcal{T}_{S}(U)={\rm Symm}(\mathcal{T}(U))$ and is called {\it the symmetric algebra generated by $U$}. It is itself a commutative unital algebra with multiplication defined as
\begin{equation*}
u \cdot v = {\rm Symm}(u \otimes v)\,.
\end{equation*}
It can be written as the span of the elementary symmetric tensors
\begin{equation*}
x_{1} \cdots x_{n}={\rm Symm}( x_{1} \otimes \cdots \otimes x_{n} ) = \frac{1}{n!} \sum_{p \in Per_{n}}  x_{p(1)} \otimes \cdots \otimes \cdots x_{p(n)}\,.
\end{equation*}

The algebra can be viewed as the free commutative algebra generated by the vector space $U$. It can also  be viewed as the quotient of the tensor algebra by the two-sided ideal $I$ generated by the set
\begin{equation*}
S= \{ (x \otimes y - y \otimes x ) \}_{x,y \in U}\,.
\end{equation*}
\end{theorem}

The symmetric algebra comes with a natural grading:
\[
\mathcal{T}_{S}(U) = \bigoplus _{t=0}^{\infty}  \bigodot^{t} U\qquad\text{where ${\rm Symm}(\bigotimes^{t} U)=\bigodot^{t} U$}.
\]

Important elements of the symmetric tensor algebra are the \emph{homogeneous} elements, namely symmetric tensor powers of a single vector.  A crucial  algebraic result for symmetric tensors is the polarization formula \cite{ryan80} that allows us to rewrite each grade as the span of homogeneous elements. This means
\begin{equation*}
\bigodot^{t} U = span (x^{t} )_{x \in U}\,.
\end{equation*}

The polarization formula is very useful and for this reason we repeat it here for completeness. The formula reads as follows
\begin{equation}\label{polar:ref}
x_{1} x_{2} \cdots x_{t} = \frac{1}{2^{t} t! } \sum_{\epsilon_{i}=\pm 1} \epsilon_{1} \cdots \epsilon_{t}(\epsilon_{1}x_{1}+\epsilon_{2}x_{2} + \cdots + \epsilon_{t}x_{t})^{t}\,.
\end{equation}

We note here that this formula is required when proving that the space of homogeneous polynomials over $U$ is isomorphic to the space of symmetric multi-linear forms of the same order. The symmetric algebra is not the only one with the above properties. The anti-symmetrization operator ${\rm ASymm}$ is defined on the tensor algebra $\mathcal{T}(U)$ as
\begin{equation}\label{asymm:def}
{\rm ASymm}( x_{1} \otimes \cdots \otimes x_{n}) := \frac{1}{n!}  \sum_{p \in Per_{n}}  sign(p) (x_{p(1)} \otimes \cdots \otimes \cdots x_{p(n)})\,,
\end{equation}
where $sign(p)$  is the sign of the permutation, see \cite{rodrigues07}. We can extend the mapping via linearity to the whole tensor algebra.

We collect in the next theorem some useful facts on the anti-symmetrization operator \cite{rodrigues07}.

\begin{theorem}\label{thm:asymm}
Let $U$ be a vector space. The anti-symmetrization operator defined by (\ref{asymm:def}) is a projection. The projection of the tensor algebra is denoted by $\mathcal{T}_{A}(U)={\rm ASymm}(\mathcal{T}(U))$ and is called the
{\it antisymmetric algebra generated by $U$}. It is itself a non-commutative unital algebra with multiplication defined as
\begin{equation*}
u \wedge v = {\rm ASymm}(u \otimes v)\,.
\end{equation*}
It can be written as the span of the elementary antisymmetric tensors
\begin{equation*}
x_{1} \wedge \cdots \wedge x_{n}={\rm ASymm}( x_{1} \otimes \cdots \otimes x_{n} ) = \frac{1}{n!}  \sum_{p \in Per_{n}} sign(p) (x_{p(1)} \otimes \cdots \otimes \cdots x_{p(n)})\,.
\end{equation*}
The algebra can be viewed as the quotient of the tensor algebra by the two-sided ideal $I$ generated by the set
\begin{equation*}
S= \{ (x \otimes y + y \otimes x ) \}_{x,y \in U}\,.
\end{equation*}
For this reason $u \wedge v = (-1)^{pq} v \wedge u$, when both $u,v$ are elementary antisymmetric tensors of non-zero grades $p$ and $q$ respectively.
\end{theorem}

The anti-symmetric algebra comes with a natural grading:
\[
\mathcal{T}_{A}(U) = \bigoplus _{t=0}^{\infty}  (\bigwedge^{t} U)\qquad \text{where $A{\rm Symm}(\bigotimes^{t} U)=\bigwedge^{t} U$}.
\]
Contrary to the symmetric algebra, the anti-symmetric algebra is finite dimensional when the vector space is finite dimensional. It also inherits the inner product from the tensor algebra when it is equipped with one.

\subsection{Review of the Hopf algebraic structure of tensor algebra}
The tensor algebra has more structure than its grading. This structure displays itself with the ability to equip it with suitable operations to turn it to a Hopf algebra. To this end we first define the co-multiplication $\Delta : \mathcal{T}(U) \to \mathcal{T}(U) \bigotimes \mathcal{T}(U)$ of the tensor algebra. The co-multiplication should be a homomorphism and because our algebra is freely generated by the elements of $U$ it is enough to  define it on the elements of $U$. For this we make the following two definitions
\begin{enumerate}
\renewcommand{\labelenumi}{(\roman{enumi})}
\item $\Delta(\mathbf{1}) = \mathbf{1} \otimes \mathbf{1}$,
\item $\Delta(x) = \mathbf{1} \otimes x + x \otimes \mathbf{1}$\qquad where $x\in U$.
\end{enumerate}
For later convenience let $T_{n}$ be the set $\{ 1,\cdots , n \}$. For a subset $S$ of $T_{n}$ with $m$ elements and an elementary tensor $u_{n}=x_{1} \otimes \cdots \otimes x_{n}$ of grade $n$ we write $u_{S}=x_{s_{1}} \otimes \cdots \otimes x_{s_{m}}$, where the elements of $S$ are arranged in increasing order. We can easily see that
\begin{equation} \label{tree:delta}
\Delta(x_{1} \otimes \cdots \otimes x_{n}) = \sum_{S \subseteq T_{n}} u_{S} \otimes u_{T_{n} \setminus S}\,.
\end{equation}
For the definition of the antipode $S : \mathcal{T}(U) \to \mathcal{T}(U)$ we similarly define (because it is a homomorphism)
\begin{enumerate}
\renewcommand{\labelenumi}{(\roman{enumi})}
\item $S(\mathbf{1}) = \mathbf{1}$,
\item $S(x) = -x$\qquad where $x\in U$.
\end{enumerate}
Finally, due to its homomorphic property the co-unit of the algebra $\epsilon : \mathcal{T}(U) \to \mathbf{K}$  can be defined on the elements of $U$ as
\begin{enumerate}
\renewcommand{\labelenumi}{(\roman{enumi})}
\item $\epsilon(\mathbf{1}) = 1$,
\item $\epsilon(x) = 0$\qquad where $x \in U$.
\end{enumerate}

Now we are ready to state a useful well-known theorem concerning the Hopf algebra structure of the tensor algebra \cite{cartier07}.

\begin{theorem}
Given a vector space $U$ with corresponding tensor algebra $\mathcal{T}(U)$, the quadruple $\langle \mathcal{T}(U) , \Delta, S , \epsilon \rangle $ is a co-commutative and co-associative Hopf algebra.
\end{theorem}

The Hopf algebra structure of the tensor algebra carries over to its symmetric algebra.

\begin{theorem}
Let $U$ be a vector space. The symmetric tensor algebra $\mathcal{T}_{S}(U)$ generated by $U$ can be equipped with the structure of a co-commutative and co-associative Hopf algebra as
\begin{enumerate}
\renewcommand{\labelenumi}{(\roman{enumi})}
\item $\Delta_{symm} (a)  =  ({\rm Symm} \otimes {\rm Symm}) \Delta(a)$\,,
\item $S_{symm}(a)  = {\rm Symm}(S(A)) $\,,
\item $\epsilon_{symm} (a)  =  \epsilon(a)$\,,
\end{enumerate}
where $a \in \mathcal{T}_{S}(U)$.
\end{theorem}

For homogeneous elements the co-multiplication takes an elegant form

\begin{equation} \label{delta:homog}
\Delta(x^{t}) = \sum_{k=0}^{t} {t \choose k} x^{t-k} \otimes x^{k}\,.
\end{equation}

The Hopf algebra structure of the tensor algebra also carries over to its antisymmetric algebra.

\begin{theorem}
Let $U$ be a vector space. The antisymmetric tensor algebra $\mathcal{T}_{A}(U)$ generated by $U$ can be equipped with the structure of a co-commutative and co-associative Hopf algebra as
\begin{enumerate}
\renewcommand{\labelenumi}{(\roman{enumi})}
\item $\Delta_{{\rm asymm}} (a)  =  ({\rm ASymm} \otimes {\rm ASymm}) \Delta(a)$,
\item $S_{{\rm asymm}}(a)  =  {\rm ASymm}(S(A)) $,
\item $\epsilon_{{\rm asymm}} (a)  =  \epsilon(a)$,
\end{enumerate}
where $a \in \mathcal{T}_{A}(U)$.
\end{theorem}

As far as the joint algebra is concerned, observe that by theorem \ref{thm:tensorhopf} it can be equipped with a Hopf algebra structure. The joint algebra and its Hopf algebraic properties will be our focus  in the next section. We conclude this section with some remarks.

\begin{remark}
The element $\Delta{a}$, where $a \in \mathcal{T}(U)$, can be identified with an element of $\mathcal{T}(U)$ as follows
\begin{equation}
\Delta(a)  =  \sum_{(a)} a_{(1)} \bigotimes a_{(2)}= \sum_{(a)} a_{(1)} \otimes a_{(2)}\,.
\end{equation}
\end{remark}

\begin{remark}
When $U$ and $V$ two vector spaces in duality. Let $a \in \mathcal{T}(U)$ and $b \in \mathcal{T}(V)$. We have the duality relations
\begin{enumerate}
\renewcommand{\labelenumi}{(\roman{enumi})}
\item $\langle {\rm Symm}(a) , b \rangle = \langle a , {\rm Symm}(b) \rangle = \langle {\rm Symm}(a) ,
{\rm Symm}(b) \rangle$,
\item $\langle {\rm ASymm}(a) , b \rangle = \langle a , {\rm ASymm}(b) \rangle=\langle {\rm ASymm}(a) ,
{\rm ASymm}(b) \rangle$.
\end{enumerate}
\end{remark}

\section{The square product}

In this section we aim to define a Laplace pairing on the joint tensor algebra. We need it in order to define a circle  product in the sense of \cite{brouder09} on the joint tensor algebra, which we call it {\it square product}.
Once our goal is completed, we will proceed to derive the ordinary circle product of the symmetric algebra defined in \cite{brouder09} as a symmetrization of the square product of the tensor algebra. A Laplace pairing on the joint tensor algebra is a symmetric bilinear form $( \cdot | \cdot )$ which should satisfy the important identity (see \cite{brouder09})
\begin{equation}\label{laplace:prop}
( a | b c) = \sum_{(a)} (a_{(1)} | b ) ( a_{(2)} | c)\qquad \text{where $a, b, c\in\mathcal{T}(U) \bigotimes \mathcal{T}(V)$}.
\end{equation}
Due to its recursive definition, in order to make it compatible with the duality already defined on the joint tensor algebra we require
\begin{enumerate}
\renewcommand{\labelenumi}{(\roman{enumi})}
\item $(\mathbf{1} | \mathbf{1} ) =1$
\item $( 1 \otimes y | x \otimes 1 )  =  ( x \otimes 1 | 1 \otimes y)   =  \langle x , y \rangle$\qquad where  $x \in U$ and $y \in V$
\item $ (x_{1} \otimes 1  | x_{2} \otimes 1 ) = ( 1 \otimes y_{1} | 1 \otimes y_{2})  =  0$\,\qquad where $x_{1},x_{2} \in U$ and $y_{1},y_{2} \in V$
\item $(x \otimes 1 | 1  \otimes 1) = (1 \otimes 1 | x  \otimes 1) = (1 \otimes 1 | 1  \otimes y)= (1 \otimes y | 1  \otimes 1)  =  0$\qquad\qquad where $x \in U$ and $y \in V $.
\end{enumerate}

One important question is whether fixing these `` initial conditions"  or the recursion, is enough to define the Laplace pairing. Moreover, we can also question the existence of the Laplace pairing. Our next theorem gives a positive answer.

\begin{theorem}\label{laplace:exists}
A Laplace pairing defined on the joint tensor algebra generated by two vector spaces $U$ and $V$ in duality compatible with the duality exists and is unique. Moreover, it is symmetric
\[
( a | b ) = (b | a) \quad\text{where $a, b\in\mathcal{T}(U)\bigotimes \mathcal{T}(V)$}.
\]
\end{theorem}

We split the proof in several steps.

\begin{lemma}
Under the hypotheses of theorem \ref{laplace:exists}, if such a Laplace pairing exists then for $m\neq n$,
\begin{equation}\label{tensor:orthog}
( x_{1} \otimes \cdots \otimes x_{n} | y_{1} \otimes \cdots \otimes y_{m} ) =0\,.
\end{equation}
\end{lemma}

\begin{proof}
We will proceed with double induction. We prove first the important fact that for $n > 0$
\begin{equation}\label{one:ortho}
( x_{1} \otimes \cdots \otimes x_{n} | \mathbf{1} ) =0\,.
\end{equation}
We argue by induction. For $n=1$ it is true because of the compatibility conditions. Suppose now that this property holds for $n$. Then, by the splitting property in equation (\ref{laplace:prop}) we have
\[
( x_{1} \otimes \cdots \otimes x_{n} \otimes x_{n+1} | \mathbf{1} ) = ( x_{1} \otimes \cdots \otimes x_{n}  | \mathbf{1} ) (  x_{n+1} | \mathbf{1} ) = 0
\]
and so the property holds for $n+1$. In a similar fashion we prove that
\[
( \mathbf{1} | y_{1} \otimes \cdots \otimes y_{m} ) =0\,.
\]
We have proved relation (\ref{tensor:orthog}) when $m=0 < n \leq 1$. Suppose now that the orthogonality holds for every pair $(n, m)$, $N\geq n > m$. We need to show the orthogonality for every pair $(n, m)$, $N+1\geq n > m$. By the induction hypothesis it is enough to prove the assertion for the pairs $(N+1, m)$ with $N+1> m$. By the properties of the Laplace pairing (\ref{laplace:prop}) and the co-multiplication (\ref{tree:delta}),
\[
( x_{1} \otimes \cdots \otimes x_{N+1} | y_{1} \otimes \cdots \otimes y_{m} ) =
\sum_{Q \subseteq T_{m}}  ( x_{1} \otimes \cdots \otimes x_{N} | y_{Q} )( x_{N+1} | y_{T_{m} \setminus Q} )\,.
\]
Observe that the  only surviving terms occur when $|Q|=N$ and $|T_{m} \setminus Q| = 1$. But then $m=N+1$ which is not possible.
\end{proof}

\begin{lemma}
Under the hypotheses of theorem \ref{laplace:exists}, if such a Laplace pairing exists
\begin{equation}
( x_{1} \otimes \cdots \otimes x_{n} | y_{1} \otimes \cdots \otimes y_{n} ) = ( x_{1} \otimes \cdots \otimes x_{n} | \sum_{p \in Per(n)} y_{p(1)} \otimes \cdots \otimes y_{p(n)})\,.
\end{equation}
\end{lemma}

\begin{proof}
Once more we argue by induction. For $n=0, 1$ the theorem obviously holds. If the theorem holds for $n$, we need prove it for $n+1$. First observe that by using (\ref{tree:delta}) we have
\[
( x_{1} \otimes \cdots \otimes x_{n+1} | y_{1} \otimes \cdots \otimes y_{n+1} ) = ( x_{1} \otimes \cdots \otimes x_{n} | \sum_{Q \subseteq T_{n+1} }  y_{Q} ) ( x_{n+1} | y_{T_{n+1} \setminus Q})\,.
\]
The RHS can be written as
\[
\sum_{k=1}^{n+1} ( x_{1} \otimes \cdots \otimes x_{n} | \sum_{ Q \subseteq ( T_{n+1}\setminus \{ k\}) } y_{Q} ) ( x_{n+1} | y_{k} )\,.
\]
Let $M_{k}$ be the totality of ``$1-1$" onto functions $g :  T_{n+1}\setminus\{ k\} \to T_{n}$. By the induction hypothesis we can rewrite the above summations as
\[
\sum_{k=1}^{n+1} ( x_{1} \otimes \cdots \otimes x_{n} \otimes x_{n+1} | \sum_{ g \in M_{k} } y_{g^{-1}(T_{n})}  \otimes y_{k} )\,.
\]
But every such $g$ is derived uniquely by restricting a $T_{n+1}$ permutation on $T_{n+1}\setminus\{ k\}$. For this reason
\[
( x_{1} \otimes \cdots \otimes x_{n} | y_{1} \otimes \cdots \otimes y_{n} ) = ( x_{1} \otimes \cdots \otimes x_{n+1} | \sum_{ g \in Per(n+1) }  y_{g^{-1}(T_{n+1})} )
\]
which proves the lemma.
\end{proof}

\begin{lemma}
We assume the hypotheses of theorem \ref{laplace:exists}. Suppose that such a Laplace pairing exists. Moreover, let $u_{i} \in \bigotimes^{n_{i}}U$ and $v_{j} \in \bigotimes^{m_{j}}V$. Then,
\begin{equation} \label{prod:prod}
( u_{1} \otimes v_{1} | u_{2} \otimes v_{2}) = ( u_{1} | u_{2}) (v_{1} | v_{2})\,.
\end{equation}
\end{lemma}

\begin{proof}
By the splitting property (\ref{laplace:prop}) and the fact that $\mathcal{T}(U)\bigotimes\mathcal{V}$ is a tensor product of Hopf algebras we have
\begin{align*}
( u_{1} \otimes v_{1} | u_{2} \otimes v_{2}) &=  \sum_{(u_{1}), (v_{1})} ( u_{1,(1)} \otimes v_{1,(1)} | u_{2}) ( u_{1,(2)} \otimes v_{1,(2)} | v_{2})\\
&= \sum_{ \substack { (u_{1}), (v_{1}) , \\ (u_{2}),  (v_{2})} }  ( u_{1,(1)}  | u_{2,(1)}) )(v_{1,(1)} | u_{2,(2)}) ( u_{1,(2)}  | v_{2,(1)} ( v_{1,(2)} | v_{2,(2)})\,.
\end{align*}
By the previous discussion $u_{1,(1)} = u_{2,(1)}=1$ and $v_{1,(2)} = v_{2,(2)}=1$. But then
$u_{1,(2)}=u_{1}$ and $u_{2,(2)}  = u_{2}$. In the same way $v_{1,(1)}=v_{1}$ and $v_{2,(1)} = v_{2}$. For this reason (\ref{prod:prod}) holds.
\end{proof}

\begin{lemma}
We assume the hypotheses of theorem \ref{laplace:exists} and the existnce of the aforementioned Laplace pairing. Suppose that $a,b \in \mathcal{T}(U) \bigotimes \mathcal{T}(V)$. We can write
$a=\sum_{i,j} u_{i} \otimes v_{j} $, with $u_{i} \in \bigotimes^{i}U , v_{j} \in \bigotimes^{j}V$
and $b=\sum_{k,l} u^{'}_{k} \otimes v^{'}_{l} $, with $u^{'}_{k} \in \bigotimes^{k}U , v^{'}_{l} \in \bigotimes^{l}V$. Then, we have
\[
(a | b )= \sum_{i,j} \sum_{k,l} i! k! \langle {\rm Symm}(u_{i}), {\rm Symm}(v^{'}_{l})\rangle \langle
{\rm Symm}(u^{'}_{k}), {\rm Symm}(v_{j})\rangle\,.
\]
\end{lemma}

\begin{proof}
This is the result of the previous discussion and of theorem ~\ref{thm:symm}.
\end{proof}

Now we have the required form of the pairing if it exists. In the course of the previous discussion we have also proven its uniqueness. In order to conclude the proof of our main theorem (theorem \ref{laplace:exists}) we have to verify that the above defined pairing has the splitting property.

{\it Proof of Theorem \ref{laplace:exists}}. Because of the previous discussion, it suffices to prove the theorem when
$a=x_{1} \otimes \cdots \otimes x_{n}$, $b=y_{1} \otimes \cdots \otimes y_{m}$ and $c=y^{'}_{1} \otimes \cdots \otimes y^{'}_{m^{'}}$. Then,
\[
(x_{1} \otimes \cdots \otimes x_{n} | y_{1} \otimes \cdots \otimes y_{m} \otimes y^{'}_{1} \otimes \cdots \otimes y^{'}_{m^{'}} )=  n! \langle x_{1}  \cdots  x_{n} , y_{1}  \cdots  y_{m} y^{'}_{1}  \cdots  y^{'}_{m^{'}} \rangle\,.
\]
Because of the symmetrization and polarization formula (\ref{polar:ref}), it is enough to prove the theorem when $a=x^{n}$, $b=y^{m}$ and $c=y^{' m^{'}}$. Then, we have
\[
(x^{n} | y^{m} \otimes y^{'}_{1} \otimes y^{' m^{'}} )=  n! \langle x , y \rangle ^{m} \langle x , y^{'} \rangle ^{m^{'}}\,.
\]
On the other hand, by (\ref{delta:homog})
\[
\sum_{(a)} (a_{(1)} | b ) (a_{(2)} | c) = \sum_{k} {n \choose k} k! (n-k)! \langle x , y \rangle ^{m} \langle x , y^{'} \rangle ^{m^{'}}\,.
\]
The symmetric nature of the pairing follows from its definition and the self-adjointness properties of the symmetrization operator. The proof of Theorem \ref{laplace:exists} is complete.

\begin{remark}
The above theorem proves the striking fact that when $U=V$ with a self-duality, the derived pairing is the same as the one derived for the symmetric algebra in \cite{brouder09}.
\end{remark}

\begin{remark}
Another striking fact is that the Laplace pairing works only with symmetric tensors. Without defining a symmetric tensor algebra we re-derived the symmetrization operator only from pure algebraic facts.
\end{remark}

Now we are in a position to present the definition of the square product.

\begin{definition}
Let $U$ and $V$ be two vector spaces in duality via the bi-linear form
$\langle \cdot , \cdot \rangle$. We define the square product for two elements
$a, b\in\mathcal{T}(U)\bigotimes\mathcal{T}(V)$ as
\begin{equation}\label{square:def}
a \square b = \sum_{(a),(b)} ( a_{(1)} | b_{(1)}) a_{(1)} \otimes b_{(2)}\,.
\end{equation}
\end{definition}

\begin{remark}
If the vector space $U$ is already in duality with itself, the joint tensor algebra is the ordinary tensor algebra with square product defined exactly as above.
\end{remark}

\begin{remark}
If the vector space $U$ is already in duality with itself, and if we restrict the product on elements of the symmetric tensor algebra, as we shall see later
\[
{\rm Symm}(a \square b )=a \circ b\qquad\text{where $a, b\in \mathcal{T}_{S}(U)$}.
\]
In the above expression  $\circ$ is the circle product defined for the specific Laplace pairing (the bilinear form).
\end{remark}

It would be convenient if the joint tensor algebra could be turned into a unital associative algebra like in the case of the circle product. We cannot require commutativity for obvious reasons. The content of the next theorem fills this gap.

\begin{theorem}
Let $U$ and $V$ be two vector spaces in duality via the bi-linear form $\langle \cdot , \cdot \rangle$. The joint tensor algebra equiped with the square product is a unital associative algebra.
\end{theorem}

\begin{proof}
Let $a,b,c \in \mathcal{T}(U,V)$. We want to prove that
\[
(a \square b) \square c= a \square ( b \square c).
\]
The LHS expression can be computed to be
\[
(a \square b) \square c= \sum_{(a),(b),(c)} (a_{(1)} | b_{(1)})  ( a_{(2,1)} b_{(2,1)}  |  c_{(1)}) a_{(2,2)} b_{(2,2)} c_{(2)}
\]
while the RHS is equal to (due to co-associativity)
\[
a \square ( b \square c)= \sum_{(a),(b),(c)} ( a_{(1)} | b_{(2,1)} c_{(2,1)}  ) (b_{(1)} | c_{(1)}) a_{(2)} b_{(2,2)} c_{(2,2)}\,.
\]
We now use the spliting property of the pairing (\ref{laplace:prop}). The above two expressions can be written as
\begin{equation*}
(a \square b) \square c= \sum_{(a),(b),(c)} (a_{(1)} | b_{(1)}) ( a_{(2,1)}   |  c_{(1,1)})  (  b_{(2,1)}  |  c_{(1,2)}) a_{(2,2)} b_{(2,2)} c_{(2)}
\end{equation*}
and
\[
a \square ( b \square c)= \sum_{(a),(b),(c)} ( a_{(1,1)} | b_{(2,1)} ) ( a_{(1,2)} | c_{(2,1)}  ) (b_{(1)} | c_{(1)}) a_{(2)} b_{(2,2)} c_{(2,2)} .
\]
respectively. Because of the co-associativity and co-commutativity of the Hopf algebras we have used, the above expressions are the same. Finally, as we can easily show, the unit of the algebra is the ordinary unit.
\end{proof}

A weak commutativity property that generalizes the one given in \cite{brouder09}, is stated in the next lemma.

\begin{lemma}
Let $a,b,c \in \mathcal{T}(V,U)$ and let the assumptions of the main theorem be fulfilled. Then, the following identity holds
\begin{equation} \label{laplace:perm}
( a \square b | c )=  ( a | b  \square c )\,.
\end{equation}
\end{lemma}

\begin{proof}
The proof is a matter of simple algebraic manipulations. We have
\begin{align*}
(a\square b | c ) &=  \sum_{(a),(b)}  ( (a_{(1)} | b_{(1)}) a_{(2)} \otimes b_{(2)} | c)\\
&= \sum_{(a),(b)} (a_{(1)} | b_{(1)}) ( a_{(2)} \otimes b_{(2)} | c)\\
&= \sum_{(a),(b),(c)} (a_{(1)}| b_{(1)}) (a_{(2)} | c_{(1)}) ( b_{(2)} | c_{(2)})\\
&= \sum_{(b),(c)} (a | b_{(1)} \otimes c_{(1)})  (b_{(2)} | c_{(2)})\\
&= ( a | b  \square c )\,.
\end{align*}
\end{proof}

Finally, we can write the square product via the ordinary product just like in \cite{brouder09}.

\begin{theorem}
Let $a,b,c \in \mathcal{T}(V,U)$ and let the assumptions of the main theorem be fulfilled. Then, the following identity holds
\begin{equation}
a \otimes b = \sum_{(a),(b)} ( S(a_{(1)}) | b_{(1)} ) a_{(2)} \square b_{(2)}\,,
\end{equation}
where $S$ is the antipode.
\end{theorem}

\begin{proof}
The proof is the same as in \cite{brouder09} for the circle product.
\end{proof}

Now we restrict our attention to the case where $U$ is equipped with an inner product and $V=U$.  We also restrict our attention to the symmetric tensor algebra. In this case the above complex calculations can be considerably simplified. We will use the notation in \cite{brouder09}.

\begin{lemma}
Let $u,v \in \mathcal{T}_{S}(V)$, then
\begin{equation}
{\rm Symm}( u \square v ) = u \circ v\,.
\end{equation}
\end{lemma}

\begin{proof}
By the definition of the square product we have
\[
{\rm Symm}( u \square v ) = \sum_{(u),(v)} {\rm Symm}(  (u_{(1)} | v_{(1)}) u_{(2)} \otimes v_{(2)}
= \sum_{(u),(v)} (u_{(1)} | v_{(1)}) {\rm Symm}(u_{(2)} \otimes v_{(2)})\,.
\]
By the definition of the circle product observe that the RHS is equal to
\begin{equation}\label{symmcirc:def}
\sum_{(u),(v)} (u_{(1)} | v_{(1)}) {\rm Symm}(u_{(2)} \otimes v_{(2)} )=u \circ v\,.
\end{equation}
\end{proof}

\section{Automorphic property of the circle and square products}

We restrict ourselves to a vector space $U$ equipped with a self-duality. We will prove that the symmetric algebra equipped with the circle product in (\ref{symmcirc:def}) is homomorphically equivalent to the symmetric algebra equipped with the ordinary product.This is not very clearly stated in \cite{brouder09} and we aim to fill this gap. We prove this property and use it as a model to extend it to the case of the square product. Let us define the mapping $\phi$ with domain and codomain equal to the symmetric tensor algebra $ \mathcal{T}_{S} (V)$ first on elementary symmetric tensors as
\begin{enumerate}
\renewcommand{\labelenumi}{(\roman{enumi})}
\item $\phi(\mathbf{1} )=\mathbf{1}$
\item $\phi(x^{n})=x \circ \cdots \circ x$\qquad\text{where $x \in U$}
\end{enumerate}
and then extend it over the symmetric algebra by linearity.

\begin{lemma}
For $u,v \in \mathcal{T}_{S}(V)$, the above mapping $\phi$ is a homomorphism as the property
\begin{equation*}
\phi(uv)=\phi(u) \circ \phi(v)
\end{equation*}
holds. Moreover, the homomorphism is injective and surjective.
\end{lemma}

\begin{proof}
The mapping is clearly linear. We shall prove the homomorphic property for elementary homogeneous tensors. For this
we use the polarization formula to prove that for $n$ vectors
$x_{1}, x_{2}, \ldots , x_{n}$  of $U$ we have
\[
\phi(x_{1}  x_{2} \cdots x_{n})= x_{1} \circ x_{2} \circ \cdots \circ x_{n}\,.
\]
Since both RHS and LHS are symmetric multilinear functions of the vectors, by the polarization formula (\ref{polar:ref})  it is enough to prove that given a linear functional $\Omega$ on the symmetric tensor algebra,
\[
\Omega ( \phi(x^{n}) ) = \Omega(x \circ x \circ \cdots \circ x )\,.
\]
But, this is true by the very definition of the homomorphism. For this reason
\[
\Omega( \phi(x_{1}  x_{2} \cdots x_{n}) )= \Omega ( x_{1} \circ x_{2} \circ \cdots \circ x_{n} )\,,
\]
and since $\Omega$ is arbitrary the above identity is valid. Observe now that
\begin{align*}
\phi(x_{1}  x_{2} \cdots x_{n} y_{1} y_{2} \cdots y_{m}) &=  x_{1} \circ x_{2} \circ \cdots \circ x_{n} \circ y_{1} \circ y_{2} \circ \cdots \circ y_{m}\\
&= \phi(x_{1}  x_{2} \cdots x_{n}) \circ \phi( y_{1} y_{2} \cdots y_{m})
\end{align*}
and by linearity
\[
\phi(uv)=\phi(u) \circ \phi(v)\qquad\text{for every $u\in \bigodot^{t} U$ and every $v \in \bigodot^{s} U$}.
\]
By linear extension we have shown that $\phi$ is homomorphic. We prove now that $\phi$ is injective. For this suppose that there exists an element $u$ of the symmetric algebra with
\[
\phi(u)= 0\,.
\]
First observe that
\[
\phi(x^{t})= x^{t}+ \text{lower grade terms}\,.
\]
By linearity for an element $u \in \bigodot^{t} U$ we have
\[
\phi(u)= u+ \text{lower grade terms}\,.
\]
For a general element $u \in \mathcal{T}_{S}(U)$ lying in the kernel of $\phi$, we have by the natural grading
\[
u= \sum_{k=0}^{m} u_{k}\qquad \text{with $u_{k}\in\bigodot^{k} U$}.
\]
Obviously $u_{k} \neq 0$. By the homomorphic property
\[
\phi(u)= \sum_{k=0}^{m} \phi(u_{k}) = u_{m} + \text{lower grade terms}=0
\]
and we conclude that $u_{k}=0$ which is a contradiction. This proves the injectivity.
It remains to show that $\phi$ is surjective. It is enough to prove the theorem for the elementary homogeneous elements . We already know that
\[
x \circ \cdots \circ x= x^{t}+ \text{lower grade terms} \iff x^{t}= x \circ \cdots \circ x + \text{lower grade terms}\,.
\]
By iterating this procedure for the low order terms we arrive at the result that
\[
x^{t}= \sum_{k=0}^{t} a_{k} x^{\circ k} = \phi ( \sum_{k=0}^{t} a_{k} x^{ k} )\,.
\]
This ends the proof of the theorem.
\end{proof}

Now we have a model to prove the same property for the square product. Namely, the fact that the tensor algebra with this square product is homomorphically equivalent to the tensor algebra with the ordinary product. Let us define the mapping $\phi$ with domain and codomain equal to the tensor algebra $ \mathcal{T} (V)$ first on elementary tensors as
\begin{enumerate}
\renewcommand{\labelenumi}{(\roman{enumi})}
\item $\phi(\mathbf{1} )=\mathbf{1}$
\item $\phi(x_{1} \otimes\cdots \otimes x_{n} ) = x_{1} \square \cdots \square x_{n}$\,,
\end{enumerate}
and then over the whole tensor algebra by linearity.

\begin{lemma}
For $u,v \in \mathcal{T} (V)$, the above mapping $\phi$ is a homomorphism as the property
\[
\phi(uv)=\phi(u) \square \phi(v)
\]
holds. Moreover, the homomorphism is injective and surjective.
\end{lemma}

The proof is analogous to the one for the circle product.

\section{The circle product for the anti-symmetric algebra}

In this section we aim to carry out the rest of the construction we mentioned in the introduction. We aim to prove that if we define
\[
u \circ v = {\rm ASymm}( u \square v)\qquad\text{where $u, v\in \mathcal{T}_{A}(U)$}
\]
then we can create yet another product on the antisymmetric algebra. As a first step we observe that
\[
x \circ y = (x | y) + x \wedge y \qquad\text{where $x, y\in V$}
\]
and
\[
x \wedge y \circ z \wedge w = (x | z) y \wedge w + (x | w) y \wedge z + (y | z) x \wedge w + (y | w) x \wedge z + x \wedge y \wedge z \wedge w\,.
\]
Obviously the new product is neither commutative , nor antisymmetric. One also observes that the pairings operate on terms of grade 0 or 1 because of anti-symmetry. We prove that the new circle product is associative.

\begin{lemma}
For elements $u, v, w\in\mathcal{T}_{A}(U)$ we have the following relation
\[
(u \circ v ) \circ w= u \circ (v \circ w)\,.
\]
\end{lemma}

\begin{proof}
It is enough to prove the theorem for elementary antisymmetric tensors. By the properties of the antisymmetric algebra outlined in theorem \ref{thm:asymm}, it is not difficult to see that
\[
(u \circ v ) \circ w= \sum_{(u),(v),(w)} (u_{(1)} | v_{(1)} ) (u_{(2,1)} \wedge v_{(2,1)}  | w_{(1)} ) u_{(2,2)} \wedge v_{(2,2)} \wedge w_{(2)}\,,
\]
while
\[
u \circ (v \circ w) = \sum_{(u),(v),(w)}  (v_{(1)} | w_{(1)} ) ( u_{(1)} | v_{(2,1)} \wedge w_{(2,1)} )
u_{(2)}  \wedge v_{(2,2)} \wedge w_{(2,2)}\,.
\]
Since the pairing is symmetric the terms inside it antisymmetric, they must be of grade at most 1 otherwise the pairing is 0. For this reason the anti-symmetrization restricted in these cases is self-adjoint because the pairing degenerates to the ordinary duality. Then we have by the co-commutativity
\[
(u \circ v ) \circ w= \sum_{(u),(v),(w)} (u_{(1)} | v_{(1)} ) (u_{(2,1)} \otimes v_{(2,1)}  | w_{(1)} ) u_{(2,2)} \wedge v_{(2,2)} \wedge w_{(2)}
\]
and
\[
u \circ (v \circ w) =
u \circ (v \circ w) = \sum_{(u),(v),(w)}  (v_{(1)} | w_{(1)} ) ( u_{(1)} | v_{(2,1)} \otimes w_{(2,1)} )
u_{(2)}  \wedge v_{(2,2)} \wedge w_{(2,2)}\,.
\]
By the definition of the square product (\ref{square:def}) and co-associativity,
\[
(u \circ v ) \circ w=  \sum_{(u),(v),(w)} (u_{(1)} \square v_{(1)} | w_{(1)} ) u_{(2)} \wedge v_{(2)} \wedge w_{(2)}
\]
while
\[
u \circ (v \circ w) =\sum_{(u),(v),(w)}  ( u_{(1)} | v_{(1)} \square w_{(1)} ) u_{(2)} \wedge v_{(2)} \wedge w_{(2)}\,.
\]
The lemma follows from (\ref{laplace:perm}), the permutation of the square product inside the Laplace pairing.
\end{proof}

It is a matter of very simple steps to prove that the new circle product is in fact an algebra product. We record in the next theorem the analogous properties of the circle product for the anti-symmetric algebra. Let us define the mapping $\phi$ with domain and codomain equal to the antisymmetric tensor algebra $ \mathcal{T}_{A} (V)$ first on elementary tensors as

\begin{enumerate}
\renewcommand{\labelenumi}{(\roman{enumi})}
\item $\phi(\mathbf{1} )=\mathbf{1}$
\item $\phi(x_{1} \wedge \cdots \wedge x_{n} ) = x_{1} \circ \cdots \circ x_{n}$\,,
\end{enumerate}
and then we extend it over the antisymmetric algebra by linearity.

\begin{theorem}
Let $U$  be a vector space in self-duality via the bi-linear form $\langle \cdot , \cdot \rangle$. The anti-symmetric tensor algebra $\mathcal{T}_{A}(U)$ equipped with the circle product is a unital associative algebra. For $u,v \in \mathcal{T}_{A}(U)$, the above mapping $\phi$ is a homomorphism as the property
\begin{equation*}
\phi(u \wedge v)=\phi(u) \circ \phi(v)
\end{equation*}
holds. Moreover, the homomorphism is injective and surjective.
\end{theorem}

The proof is analogous to the one for the circle product on symmetric tensors. The reader can see that this circle product is indeed the product presented in \cite{rota94}.

\section{Conclusion}

In this paper we have proved that the circle product on the symmetric tensor algebra, introduced in \cite{brouder09}, can be derived as a restriction of a more general product on the tensor algebra. This general product shares many commonalities like the circle product. It is defined via a more general Laplace pairing. For symmetric tensors this Laplace pairing makes the circle multiplication by a symmetric tensor self-adjoint. For general tensors the adjoint map of left square multiplication is the left square multiplication in the opposite algebra. Finally, we have proved the important fact that the new algebraic structures are homomorphically equivalent with the traditional structures. As a by-product we have shown that the symmetric tensor algebra arises naturally via algebraic considerations, namely via the Laplace pairing properties. We have also extended our results to the anti-symmetric algebra to re-derive the product defined by Rota and Stein \cite{rota94}. This step gives us a unified picture of the circle products defined in the literature as specializations of a general structure.



\bibliographystyle{amsplain}

\begin{thebibliography}{99}


\bibitem[1]{brouder09}
C. Brouder, {\em Quantum field theory meets Hopf algebra}, Math. Nachr. 282 (2009), pp. 1664--1690.

\bibitem[2]{cartier07}
P. Cartier, {\itshape A primer of Hopf algebras}, Frontiers in number theory, physics, and geometry. II, Springer, Berlin, 2007, pp. 537--615.

\bibitem[3]{fauser01}
B. Fauser, {\em On the Hopf algebraic origin of Wick normal ordering}, J. Phys. A 34 (2001),
pp. 105--115.

\bibitem[4]{fauser03}
B. Fauser, {\em Quantum Clifford Hopf gebra for quantum field theory}, Adv. Appl. Clifford Algebras 13 (2003),
pp. 115--125.

\bibitem[5]{kassel95}
C. Kassel, {\itshape Quantum groups}, Graduate Texts in Mathematics, 155, Springer-Verlag, New York, 1995.

\bibitem[6]{majid02}
S. Majid, {\itshape A quantum groups primer}, London Mathematical Society Lecture Note Series, 292,
Cambridge University Press, Cambridge, 2002.

\bibitem[7]{podles98}
P. Podle\'s and E. ~M\"uller, {\em Introduction to quantum groups}, Rev. Math. Phys. 10 (1998), pp. 511--551.

\bibitem[8]{rodrigues07}
W.A. Rodrigues, Jr. and E. Capelas de Oliveira, {\itshape The many faces of Maxwell, Dirac and Einstein equations},
A Clifford bundle approach, Lecture Notes in Physics, Vol. 722, Springer 2007.

\bibitem[9]{rota94}
G.-C. Rota and J.A. ~Stein, {\em Plethystic Hopf algebras},  Proc. Nat. Acad. Sci. U.S.A. 91 (1994), pp. 13057--13061.

\bibitem[10]{ryan80}
R. A. Ryan, {\em Applications of topological tensor products to infinite dimensional holomorphy}, Ph. D. Thesis,
Trinity College Dublin, 1980.

\end{thebibliography}

%
%

\end{document}